\documentclass{amsart}

\usepackage{amssymb,amsmath}

\def\bc{\begin{center}}
\def\ec{\end{center}}
\def\be{\begin{equation}}
\def\ee{\end{equation}}
\def\ben{\begin{enumerate}}
\def\een{\end{enumerate}}
\def\bfg{\begin{figure}}
\def\efg{\end{figure}}
\def\bq{\begin{quote}}
\def\eq{\end{quote}}
\def\bd{\begin{description}}
\def\ed{\end{description}}
\def\this{i.\ e.\ } %that is

\def\h{\hbar}
\def\p{\partial}
\def\w{\wedge}

\def\dim{\operatorname{dim}}

\def\det{\operatorname{det}}

\newcommand{\CC}{{\Bbb C}}

\newcommand{\QQ}{{\Bbb Q}}

\newcommand{\lan}{\langle}
\newcommand{\ran}{\rangle}
\newcommand{\gel}{\varepsilon}

\newcommand{\gl}{\lambda}

\newcommand{\gr}{\rho}

\newcommand{\gt}{\tau}

\newcommand{\M}{\overline{\mathcal M}}
\newcommand{\ct}{\operatorname{ct}}
\newcommand{\ev}{\operatorname{ev}}
\renewcommand{\QQ}{\mathbf Q}
\renewcommand{\CC}{\mathbf C}
\renewcommand{\t}{\mathbf t}
\newcommand{\T}{\mathbf T}
\newcommand{\V}{\mathbf V}
\newcommand{\Dl}{\mathbf D}
\renewcommand{\a}{\alpha}
\renewcommand{\b}{\beta}
\renewcommand{\c}{\gamma}
\renewcommand{\d}{\delta}
\newcommand{\x}{\chi}

\begin{document}

\title{Semisimple Frobenius structures \\ at higher genus}
\author{Alexander B. Givental} 
\address{UC Berkeley 
         and Caltech}  

\dedicatory{To the memory of Tom Wolff}
\thanks{Research partially supported by NSF Grant DMS-0072658} 
\date{August 8, 2000, revised May 4, 2001}

%\begin{abstract}
%In the context of equivariant Gromov-Witten theory of
%tori actions with isolated fixed points we compute
%genus $g\geq 2$ Gromov-Witten potentials and their
%generalizations with gravitational descendents.
%Both formulas, with and without descendents, are
%stated in a form applicable to any semisimple Frobenius
%structure and therefore can be considered as 
%definitions in the axiomatic context of Frobenius manifolds. 
%In (non-equivariant) Gromov-Witten theory, they become 
%conjectures expressing higher genus GW-invariants in terms
%of genus $0$ GW-invariants of symplectic manifolds with 
%genericly semisimple quantum cup-product.
%\end{abstract}

\maketitle

\section*{Introduction}

{\em The genus $g$ GW-potential} of a compact symplectic manifold $X$
is a generating function for genus $g$ Gromov-Witten invariants. 
It is a formal function 
\[ F^g_X(t):=\sum_{n=0}^{\infty}\sum_{d\in H_2(X)}\frac{q^d}{n!}  
\int_{[X_{g,n,d}]}\ev^*_1(t)\w ... \w \ev_n^*(t) ,\]
on the cohomology space $H^*(X,\QQ \{\{q\}\})$ 
over a suitable Novikov ring $\QQ \{\{q\}\}$. The coefficients are defined 
by integration over virtual fundamental cycles in the moduli spaces of 
degree $d$ genus $g$ stable pseudo-holomorphic curves with $n$ marked points. 
The cohomology classes $\ev_i^*(t)$ are pull-backs from $X$ by the evaluation 
maps at the marked points. 

One may use the natural {\em contraction maps} $\ct: X_{g,n,d}\to \M_{g,n}$
to the Deligne -- Mumford moduli spaces of marked Riemann surfaces in order to
define more general potentials by integration over inverse images of boundary
strata or of  any other cycles.

The potentials $F^g_X$ and their generalizations are expected to obey 
some universal constraints, yet unknown explicitly 
(see however \cite{BP,Ge,Ge2,KM}), 
but encoded implicitly in
the topology of the Deligne-Mumford spaces $\M_{g,n}$. In a sense, the 
implicit constraints, to be considered as axioms of $2$-dimensional 
Topological Field Theory, are the subject of our study in this paper.         

\medskip

In this paper, we will compute genus $g\geq 2$ Gromov-Witten invariants and their generalizations with gravitational descendents in the context of equivariant Gromov-Witten theory of tori actions with isolated fixed points. 
Both formulas, with and without descendents, are
stated in a form applicable to the axiomatic version of
genus $0$ Gromov -- Witten theory, namely --- to semisimple
Frobenius structures.  Therefore the formulas can be considered 
as definitions extending the genus $0$ theory to higher genus in a way
consistent --- conjecturally --- with the implicit axioms mentioned above.
In (non-equivariant) Gromov-Witten theory, the formulas become 
conjectures expressing higher genus GW-invariants in terms
of genus $0$ GW-invariants of symplectic manifolds with 
genericly semisimple quantum cup-product.

\section{Definitions and examples}

\subsection{ Frobenius structures.}

The axiomatic structure of $2$D TFT is understood well in genus $0$ due to 
R. Dijkgraaf -- E. Witten \cite{DW}, B. Dubrovin \cite{Db} and many others
(see \cite{Mn}) as the theory of Frobenius manifolds. By definition, 
a {\em Frobenius structure} on a manifold $H$ consists of:
 
\noindent (i) a flat pseudo-Riemannian metric $(\cdot ,\cdot)$,
    
\noindent (ii) a function $F$ whose $3$-rd covariant derivatives 
$F_{abc}$ are structure constants $(a\bullet b,c)$ of a Frobenius algebra
structure, i.e. associative commutative multiplication
$\bullet $ satisfying $(a\bullet b, c)=(a,b\bullet c)$, on the tangent spaces
$T_tH$ which depends smoothly on $t$; 

\noindent (iii) the vector field of unities ${\bf 1}$ of the $\bullet$-product
which has to be covariantly constant and preserve the multiplication 
and the metric.  

\medskip

{\em Example 1.} The genus $0$ GW-potential $F=F^{0}_X$ defines a Frobenius
structure on the super-space  $H=H^*(X,\QQ)$ 
\footnote{Formally speaking $F^{0}_X$ defines a Frobenius structure over $\QQ \{\{q\}\}$.
However, due to the divisor equation, $3$-rd rderivatives of $F^{0}_X$ make sense at $q=1$
as formal Fourier series along $H^2(X,\QQ)$ and define a Frobenius structure over $\QQ$. 
We refer to \cite{Au,G1,Mn} for discussions of these standard subtleties.} 
In this example, 
the metric and the unit vector field are translation-invariant and defined 
by the Poincare intersection pairing and by the cohomology class $1$ 
respectively.

\medskip

{\em Example 2.} Let $f(x,t)$, $t\in H$, 
be a miniversal deformation (with respect to 
the right equivalence) of the germ $f(\cdot, 0): (\CC^m,0)\to (\CC,0)$ of 
a holomorphic function at an isolated critical point.  
Then the tangent spaces $T_tH$ are canonically identified with the
algebras $Q_t:=\CC\{x\}/(f_x)$ of functions on the critical schemes 
$\operatorname{crit} f(\cdot,t)$ and thus carry a natural multiplication
$\bullet$ with unity $1$. 
Let $\Omega$ be a holomorphic volume form on $\CC^m$ possibly
depending on $t$. The multiplication $\bullet $ is Frobenius with respect to
the residue pairing
\[ (\phi ,\psi):= \frac{1}{(2\pi i)^m} \oint_{|f_{x_1}|=\gel_1} ... 
\oint _{|f_{x_m}|=\gel_m} \frac{\phi (x) \psi (x)\ \Omega}
{f_{x_1}...f_{x_m}},\]
which is known to be non-degenerate on $Q_t$ (see \cite{GH}).
According to the theory \cite{Sa} of {\em primitive volume forms} 
there exists a choice of $\Omega$ such that the corresponding residue metric
is flat and constitutes, together with the multiplication $\bullet$, 
a Frobenius structure on $H$ (see also \cite{Ba} for a new approach).

\medskip

Frobenius manifolds of Examples $1$ and $2$ come equipped with
one more ingredient --- the {\em Euler vector field} $E$ such that  
$\bullet ,\ {\bf 1}$ and $(\cdot,\cdot )$ are eigenvectors of the Lie 
derivative $L_E$ with the eigenvalues $0$, $-1$ and $2-D$ respectively.
Such Frobenius structures are called {\em conformal}, and $D$ is called their
{\em dimension}. In the Example $1$,
$D$ coincides with the complex dimension of the target manifold $X$, and the
grading imposed by $E$ originates from grading in cohomology. In Example $2$,
the Euler vector $E(t)$ is given by the class of the function $f(\cdot,t)$
in the algebra $Q_t$, and $D=1-2/h$ where $h$ is the so called {\em Coxeter
number} of the singularity \cite{AVG}. Frobenius manifolds in the next 
example fall out of the conformal class.

\medskip

{\em Example 3.} Let the K{\"a}hler manifold $X$ be endowed with 
a Hamiltonian Killing action of a compact group $T$. Then one
can introduce {\em equivariant GW-invariants} \cite{G1} using $T$-equivariant
cohomology and intersection theory in the moduli spaces $X_{g,n,d}$. The 
genus $0$ equivariant GW-invariants define on $H:=H^*_T(X,\QQ)$ 
the structure of a Frobenius manifold {\em over the ground ring} 
$H^*(BT,\QQ)$, the coefficient ring of the equivariant cohomology theory.
On the other hand, grading in equivariant cohomology imposes homogeneity
constraints on GW-potentials so that $(\cdot,\cdot),\ {\bf 1}$ and $\bullet$ do
have degrees $2-\dim X$, $-1$ and $0$ with respect to a suitable Euler
vector field $E$. Yet the Frobenius structure is not conformal since
elements of the ground ring may have non-zero degrees and therefore
$L_E$ is a differentiation only over $\QQ$ instead of the ground ring
of the Frobenius structure.

\medskip

A Frobenius manifold is called {\em semisimple} if the algebras 
$(T_tH,\bullet)$ are semisimple at generic $t$. Frobenius structures of
Example $1$ are semisimple for, say, projective spaces and flag manifolds,
and are not semisimple for Calabi-Yau manifolds. Let us assume now on
that the group $T$ in Example $3$ is a torus acting on $X$ with isolated fixed
points only. Then the cup-product in the equivariant cohomology $H^*_T(X,\QQ)$
is genericly semisimple, resulting in the corresponding Frobenius structure 
being semisimple too. All Frobenius manifolds of Example $2$ are semisimple.

\subsection{  The formula.}

Our expression for the higher genus potentials $F^g$ of a semisimple
Frobenius manifold $H$ has the form
\begin{equation} \label{1}
 \begin{array}{ll}  
e^{ \sum_{g\geq 2}\h^{g-1}F^g(t)} \ =\ & \ \\ 
=\ [ \ e^{\frac{ \h}{2}\sum_{k,l=0}^{\infty} \sum_{i,j} V^{ij}_{kl}
\Delta_i^{1/2}\Delta_j^{1/2}\p_{Q^i_k}\p_{Q^j_l} }  \  
\prod_{j}\ \tau (\h \Delta_j;\ Q^j_0,Q^j_1,...)\ ]
\ _{Q^i_k=T^{i}_k}, & \ \end{array} 
\end{equation}
where $ V^{ij}_{kl}, \Delta_j, T^i_k$
are certain functions of $t\in H$ defined at semisimple points, 
$i,j=1,...,\dim H$, $k,l=0,1,2,...$, 
and $\tau $ is the following {\em Kontsevich -- Witten tau-function}.

Let $c^{(1)},...,c^{(n)}$ denote the $1$-st Chern classes 
of the {\em universal cotangent lines} over the Deligne -- Mumford spaces 
$\M_{g,n}$, i.e. line bundles formed by cotangent lines to the
curves at the marked points. We put $Q(c)=Q_0+Q_1c+Q_2c^2+...$ 
where $Q_i$ are formal variables,
introduce the {genus $g$ descendent potential} of $X=\operatorname{pt}$ 
\[ {\mathcal F}^g_{\operatorname{pt}}(Q)=\sum_{n=0}^{\infty} 
\frac{1}{n!}\int_{\M_{g,n}} Q(c^{(1)})\w ...\w Q(c^{(n)}) \]
and define
\begin{equation} \label{2} \tau (\h ; Q)= \exp \{ \sum_{g=0}^{\infty} \h^{g-1} 
{\mathcal F}^g_{\operatorname{pt}} (Q) \} .\end{equation}
As it was proved by M. Kontsevich 
\cite{Kn}, $\tau (Q)$ provides an asymptotic expansion of the matrix Airy
function and (modulo some re-notation) coincides, as it was conjectured by E. Witten
\cite{W}, with the tau-function of the KdV-hierarchy satisfying the string equation.

In order to define the functions $V^{ij}_{kl}, \Delta_i, T_k^j$ 
we have to review the structural theory of semisimple Frobenius manifolds 
\cite{Db,G2,Mn}.

\subsection{Canonical coordinates, Hessians and stationary phase asymptotics.}

Given a germ of a Frobenius manifold, we introduce coordinates $\{ t^{\a} \}$
flat with respect to the metric $(\cdot,\cdot)$, denote $\{ \phi_{\a} \}$ 
the corresponding frame in the tangent bundle, put $g_{\a\b}:=
(\phi_{\a},\phi_{\b})$ and $(g^{\a\b}):=(g_{\a\b})^{-1}$.

The associativity constraint of the $\bullet$-product is expressed by 
the {\em WDVV-identity} for the genus $0$ potential (we use the summation 
convention if possible): 
\[ F_{\a\b\mu}g^{\mu\nu}F_{\nu\c\d}=F_{\a\c\mu}g^{\mu\nu}F_{\nu\b\d}=
F_{\a\d\mu}g^{\mu\nu}F_{\nu\b\c}.\]
It can also be interpreted as commutativity of the following
connection operators $\nabla_{\a}(z):= z\p_{\a}+\phi_{\a} \bullet$ on $TH$
and respectively --- the compatibility property of the following
linear PDE system on $T^*H$ for any value of the parameter $z\neq 0$: 
\begin{equation} \label{3}
  z\p_{\a}S_{\b}=F_{\a\b\mu}g^{\mu\nu}S_{\nu}.\end{equation}

The system plays an important role in the theory of Frobenius
structures, and we ought to start with some remarks about its solutions.  
A fundamental solution to (\ref{3}) can be found in
the form of a power $z^{-1}$-series $1+z^{-1}S_1+z^{-2}S_2+...$ satisfying 
the unitary condition $ S^*(-1/z) S (1/z) =1$ (here ``$\ ^*$'' means ``adjoint
relative to $(\cdot,\cdot)$").
Such a solution $S$ is unique up to right multiplication by a constant
matrix $1+O(z^{-1})$ satisfying the unitary condition. A choice of such 
a solution is the starting point in Dubrovin's construction \cite{Db} of genus $0$
gravitational descendents of Frobenius structures. We will review this construction in 
the section $3$ and respectively will make use of such a fundamental solution in
our description (\ref{22}) of higher genus {\em descendent} potential. 
However, both higher genus formulas (with or without descendents) require another,
asymptorical form of solution to the same system (\ref{3}) which can be constructed
for a {\em semisimple} Frobenius structure as follows.     

Let us assume that the Frobenius manifold is semisimple. In a neighborhood
of a semisimple point one introduces {\em canonical coordinates}
$\{ u^i(t) \}$ (see \cite{Db}).
They are characterized uniquely up to reordering and additive constants by
the property of $\p_i:=\p/\p u^i$ to form the basis of canonical idempotents
of the $\bullet$-product on $T_tH$. 
The flat metric $(\cdot,\cdot)$ is diagonal in canonical coordinates
and is therefore determined by the non-vanishing functions $(\p_i,\p_i)$.
We put $\Delta_i:= 1/(\p_i,\p_i)$. In singularity theory, $u_i$ are
critical values of the Morse functions $f(\cdot,t)$ at the critical
points, and $\Delta_i$ are the Hessians at these points computed in 
$\Omega$-unimodular coordinate systems.

Let $U$ denote the diagonal matrix of canonical coordinates 
$\operatorname{diag}(u_1,...,u_N)$, and $\Psi $ denote the transition
matrix between the flat and normalized canonical bases:   
$\Delta_i ^{-1/2} du ^i=\sum_{\b}\Psi_{\b}^i dt^{\b}$.
In particular, $\sum_i\Psi_{\a}^i\Psi_{\b}^i=g_{\a\b},\ 
\Psi_{\mu}^ig^{\mu\nu}\Psi_{\nu}^j=\d_{ij}$.

\medskip

{\bf Proposition} (see \cite{Db,G2}). 

{\em (a) Near a semisimple point the system (\ref{3}) has a fundamental
solution in the form of the matrix series :
\begin{equation} \label{4}
S = \Psi (R_0+z R_1+z^2R_2+...)\exp {U/z} \end{equation}
where $R_k=(R_{k})_i^j$ are matrix-functions of $u$, and $R_0=1$.

(b) The series solution $S$ can be chosen to satisfy the unitary
condition 
\begin{equation}  
(1+zR_1+z^2R_2+...)(1-zR_1^t+z^2R_2^t-...)=1 \end{equation}

(c) The series $R=1+zR_1+z^2R_2+...$ in the solution $S$ satisfying 
the unitary condition is unique up to right multiplication by 
unitary diagonal matrices $\exp (a_1z+a_2z^3+a_3z^5+ ...)$  where
$a_{k}=\operatorname{diag}(a_{k}^1,...,a_{k}^N)$ are  
constant.

(d) In the case of conformal Frobenius structures the
series $R$ in a fundamental solution $S$ can be chosen homogeneous, and
 such $R$ is unique and possesses the unitary property automatically.}

\medskip

{\em Proof.} A proof of (d) and (a) is given in \cite{Db} and 
\cite{G2}. We will remind below some details from \cite{G2}
in order to justify the additions (b) and (c) needed here.

Substitution of $S=\Psi (1+... )\exp(U/z)$ into (\ref{3}) yields
a chain of equations $(d+W\w) R_{k-1}=[dU,R_{k}]$, where 
$W=\Psi^{-1}d\Psi=[dU,R_1]$, to be solved inductively starting with
$R_0=1$. First, off-diagonal entries of $R_k$ are expressed
algebraically via $R_{k-1}$, then the diagonal terms of $R_k$ are
found by integration from the next equation using the fact that
$[dU,R_{k+1}]$ has zero diagonal entries. Compatibility conditions needed in
this procedure are verified in \cite{G2}. 

In order to prove (b), let us introduce a temporary notation
$P_k=R_kR_0^t-R_{k-1}R^t_1+...(-1) ^kR_0R^t_k$ for the $z^k$-term in
$R(z)R^t(-z)=1+P_1z+P_2z^2+...$. A short elementary computation shows
that $[dU, P_k]=d P_{k-1}+[W, P_{k-1}]$. Assuming that $P_{k-1}=0$ 
(or $1$ for $k=0$), we conclude that off-diagonal entries of $P_k$
vanish. This already implies $P_k=0$ for odd $k$ since such $P_k$ are
obviously anti-symmetric. Now, taking in account that $P_k$ is diagonal and 
$W=\Psi^{-1}d\Psi$ is anti-symmetric, we conclude from
the next equation $dP_k+[W,P_k]=[dU,P_{k+1}]$ that the 
diagonal entries of $P_k$ are constant. For even $k$ we have 
$P_k=R_k+R_k^t+...$ and thus a unique choice of integration 
constants in the above procedure for finding $R_k$ will make $P_k$ vanish.
 
Yet the integration constants for diagonal entries of $R_{2k-1}$ 
are totally ambiguous, and it is immediate to see, by induction 
on $k$, that this ambiguity is correctly accounted by the multiplication
$R\mapsto R\exp (a_{k}z^{2k-1})$ described in (c).   
   
In the conformal case, let $E=\sum u ^i\p_i$ denote the Euler
field. The Euler formula $R_k=-(i_EdR_k)/k$ shows how to recover
diagonal entries of $R_k$ via their differentials by an algebraic
procedure. This implies existence of a homogeneous solution $R$. 
Finally, the homogeneity condition leaves no freedom in the choice 
of the integration constants, but it also guarantees that the constant 
diagonal entries in $P_{2k}$ are zeroes. This proves (d). 

\medskip

Let $S(z)$ be the unitary fundamental solution to
(\ref{3}) singled-out in the proposition. 
We introduce a new matrix-function 
\[ [V^{ij}(z,w)] :=(z+w)^{-1} [S_{\mu}^i(z)]^t [g^{\mu\nu}]
[S_{\nu}^j(w)] .\]
It expands as $V^{ij}(z,w)=$
\begin{equation} \label{6} 
\frac{e ^{u ^i/z+u^j/w}}{z+w}\sum_{s} R^i_s(z)R^j_s(w) 
 =: e ^{u ^i/z+u^j/w} 
(\frac{\d^{ij}}{z+w} + \sum_{k,l=0}^{\infty}(-1)^{k+l} V^{ij}_{kl}z^kw^l)
.\end{equation}
This defines $V_{kl}^{ij}$ as functions on the Frobenius manifold 
in a neighborhood of a semisimple point.
 
Next, in the semisimple Frobenius algebras $(T_tH, \bullet)$ we have~:
\[ 1=\sum \d^{\mu} \phi_{\mu} = \sum \p_j= \sum \Delta_j ^{-1/2}
(\Delta_j ^{1/2} \p_j) .\]
We expand "the first row" of $S(z)$
\begin{equation} \label{7}
 \sum \d^{\mu} S_{\mu} ^i(z) =
(\sum_j \Delta_j ^{-1/2} R_j^i(z)) e^{u^i/z}=: 
[1- \sum_{k=0}^{\infty} T^i_{k}(-z)^{k-1} ] \frac{e^{u^i/z}}{\sqrt{\Delta_i}}. 
\end{equation}
This defines $T_k^i$. In particular, $T_0^i=T_1^i=0$.

\medskip

Perhaps, the nature of these formulas, the asymptotical solution $S(z)$ and the relationship between the two forms of fundamental solutions to (\ref{3}) will become more transparent after the following two examples. 

\medskip

{\em Example 4.} In singularity theory, a fundamental solution
matrix to the equation (\ref{3}) is given by
complex oscillating integrals of suitable $m$-forms over suitable $m$-cycles~:
\[ S_{\mu}^{i}=\int_{\Gamma ^{i}\subset \CC^m} e^{f(x,t)/z}\phi _{\mu}(x,t)
\Omega \ .\]
The cycles $\Gamma^i$ can be constructed as in Morse theory for
the function $\operatorname{Re} \{ f(\cdot,t)/z \} $ and thus correspond 
to critical points $x^i(t)$ of the function $f(\cdot,t)$.
To construct $1/z$-expansion of the integrals, one first expands the
integrals over the levels $f=f_{crit}-\tau $ near $\tau =\infty$
and then describes the oscillating integrals via the Laplace transform.
Say, for weighted-homeogeneous singularities 
\[ \int_{\c_i} \phi_{\mu}(x,t) \frac{\Omega}{df (x,t)} = 
\tau^{d_{\mu}}(\sum A_{k,\mu}^i(t) \tau^{-k}) ,\]
where $d_{\mu} $ is the weight of the form $\phi_{\mu} \Omega/df$.
Respectively,
\[ S_{\mu}^i = \sum A_{k,\mu}^i \int_{0}^{\infty} e^{-\tau/z} \tau^{d_{\mu}-k} \ d\tau 
= z^{d_{\mu+1}} \sum A_{k,\mu}^i(t) \Gamma (d_{\mu}+1-k) z^{-k} .\] 

Alternatively, one arrives to the expansion (\ref{4}) via the stationary phase
asymptotics of the oscillating integrals near non-degenerate critical points of Morse
functions $f(\cdot,t)$:
\[ \int_{\Gamma^i}e^{f(x,t)/z}\phi_{\mu}(x,t)\ \Omega \sim 
e^{u^i/z}(\frac{\phi_{\mu}(x^i,t)}{\sqrt{\Delta_i}} + ... ) \]
where $u^i=f(x^i,t)$ is the critical value and $\Delta_i$ is the
$\Omega$-Hessian at the critical point.
In particular (\ref{7}) is the stationary phase expansion  
\[ \int_{\Gamma^i}e^{f/z}\Omega \sim \frac{e^{u^i/z}}{\sqrt{\Delta_i}}
[1+T_2^iz-T_3^iz^2+... ] \]

{\em Example 5.} In Gromov-Witten theory, a $1/z$-series
solution  to (\ref{3}) satisfying the unitary condition is given by the following 
matrix of gravitational {\em descendents}:
\begin{equation}
\label{8}
 \lan \phi_{\b}, \frac{\phi_{\c}}{z-c}\ran:=
\sum_{n,d}\frac{q^d}{n!}\int_{[X_{0,2+n,d}]} 
\ev_{0}^*(\phi_{\b})\w \ev_1^*(t)\w...\w\ev_n^*(t)\w
\frac{\ev_{n+1}^*(\phi_{\c})}{z-c^{(n+1)}} \ .\end{equation}
By definition, the constant $g_{\a\b}$ is taken on the role of the
ill-defined term with $d=0, n=0$. 
This solution is related to the {\em two-point descendent} 
\begin{equation} \label{9}
 \lan \frac{\phi_{\a}}{z-c},\frac{\phi_{\b}}{w-c}\ran :=
\sum_{n,d}\frac{q^d}{n!}\int_{[X_{0,2+n,d}]}
\frac{\ev_0^*(\phi_{\a})}{z-c^{(0)}}\w \ev_1^*(t)\w...\w\ev_n^*(t)\w
\frac{\ev_{n+1}^*(\phi_{\b})}{w-c^{(n+1)}}.  \end{equation}
in the same way as $S(z)$ is related to $V(z,w)$:
\begin{equation} \label{10}
\lan \frac{\phi_{\a}}{z-c},\frac{\phi_{\b}}{w-c}\ran =\sum_{\mu\nu}
\lan \frac{\phi_{\a}}{z-c},\phi_{\mu}\ran g^{\mu\nu} \lan \phi_{\nu}, \frac{\phi_{\b}}{z-c}\ran
. \end{equation}
According to the mirror conjecture
\cite{G3,G2} the descendents (\ref{8}) can be identified with oscillating integrals
of the mirror partner. When this is the case the values of $\Delta^i$ and $T_k^i$ can 
be extracted from the stationary phase asymptotics of 
the integrals. 
In Section $2$ (see also \cite{G2}), we will find that in the equivariant setting of Example $3$
when the fixed points of the toris action are isolated and respectively the classical 
equivariant cohomology algebra of the target space is semisimple, the series
$S(z)$ and $V(z,w)$ essentially coincide with the descendents (\ref{8}) and (\ref{9}).  

\medskip

{\bf Conjecture $1$.} {\em With the notations (\ref{2},\ref{6},\ref{7})
in force, the formula (\ref{1}) represents higher genus GW-invariants
of compact symplectic manifolds with generically semisimple quantum 
cup-product.}

\medskip

The main reason to believe in the conjectural
formula (\ref{1}) is the theorem in the next section and the
miraculous coincidences which occur in the proof. We would like to mention here one
more bit of evidence in its favor. Namely, in the case of conformal
Frobenius manifolds of dimension $N=2$ (``two primaries'') 
our formula yields, after some computation, the following genus $2$ potential:
\[ \frac{d(3d-1)(d-1)^2(3d-5)(d-2)}{2880}\ \frac{\Delta_{-}}{(u^{+}-u^{-})^{3}}
, \]
where $d$ is the conformal dimension, and $u^{\pm}$ are the canonical
coordinates. This answer coincides with the result found in \cite{EGX}.

Note that the potential vanishes in the case $d=1/3$ 
corresponding to the singularity of type $A_2$. 
This fact agrees with the general conjecture that our formula (\ref{1}), 
when applied to the Frobenius
structures on the miniversal deformations of isolated critical points, 
should give rise to higher genus potentials which extend analytically 
through the bifurcation set 
(and therefore must vanish for $A,D,E$-singularities). 
For no apparent reason, the above formula is symmetric about $d=1$ 
(corresponding to GW-invariants of $\CC P^1$). It would be interesting 
to find out what is behind this symmetry.

\section{Computation in equivariant GW-theory}

In this section, we formulate and prove Theorem $1$ confirming the
conjectural formula (\ref{1}) in the case of equivariant Gromov -- Witten invariants
of Hamiltonian tori actions with isolated fixed points. Roughly speaking,
we will compute the GW-invariants using fixed point localization and
will see how the formula (\ref{1}) emerges from the combinatorial
formalism of summation over graphs. We will first discard those
factors in localization formulas which are due to the so called
Hodge intersection numbers. This will lead us to a (wrong!) 
higher genus potential formula based on a matrix series $R(z)$ 
corresponding to some reference choice $a_k^i=0$ of the integration 
constants of the part (c) of Proposition. Then we will point out
a new choice of the integration constants $a_k^i$  which compensates
for the effect of the Hodge integrals and yields a right formula
for the higher genus potential. 

\subsection{  Localization and materialization.}
 
Let the torus $T$ act on $X$ with isolated 
fixed points only. Fixed points of the induced action of $T$ on the
moduli spaces $X_{g,n,d}$ can be described as curves formed by {\em legs}
--- $1$-dimensional orbits of $T_{\CC}$ in $X$ or their multiple covers, ---
which are connected at {\em joints} --- nodes or DM-stable curves mapped 
to fixed points $X^T$. Due to multiplicative properties of the Euler classes
contributions of fixed points into localization formulas 
essentially factors into contributions of legs and joints \cite{Kn1,GP}.
(We are assuming for simplicity that the $1$-dimensional
orbits are also isolated. Beyond this assumption, our arguments remain valid 
but the leg contributions are to be found by integration over suitable orbit
spaces.)
Contributions of fixed point submanifolds can be arranged as the sum 
over strata in Deligne -- Mumford spaces in accordance with 
images of the submanifolds under the contraction map 
$\ct: X_{g,n,d}^T\to \M_{g,0}$. It is
convenient to name some elements of $T$-invariant curves depending on their
fate under the contraction map. We call {\em vertices} those joints of
$T$-invariant curves in $X$ which contract to irreducible components of 
DM-stable $(g,0)$-curves. The genus $0$ trees of legs and joints which 
contract to (self-)intersection points of these components are called 
{\em edges}. The trees which contract to non-singular points are called
{\em tails}. 

Thinking of a $T$-invariant curve (may be disconnected)
as a collection of vertices (DM-stable curves mapped to the fixed points 
$X^T$) with arbitrary number of tails attached and connected somehow by the 
edges, we arrive at the fixed point expression for the higher genus potential
with the standard combinatorics (\ref{1}) of Wick's formula. 
Contributions of vertices will be expressible via intersection numbers 
(\ref{2}) in Deligne -- Mumford spaces, while the edge factors and tail factors
should be extracted from genus $0$ GW-invariants of $X$.

A key point is that {\em the genus $0$ data needed in the localization formulas
can be written in abstract terms of semisimple Frobenius structures, 
and vice versa.} For example, in the GW-theory of $X$, the sum $\sum u^i$ 
of canonical coordinates enumerates elliptic curves with a
{\em fixed complex structure}. Expressing the GW-invariant
via the sum over fixed point components we can single out the sub-sum
where {\em the} elliptic joint of the curve is mapped to the $i$-th
fixed point in $X$. It turns out \cite{G1} that the sub-sum equals $u_i$.
Another example: let $\{\phi_{\a}\}$ be the basis of $\d$-functions at the 
fixed points in localization of $H^*_T(X)$, so that 
$g^{\a\b}=e_{\a}\d_{\a\b}$ where $\sum e_{\a}\phi_{\a}$ is the 
equivariant Euler class of $TX$. 
In the fixed point sum for $F^{(0)}_{\a\a\a} e_{\a}^{3/2}$
(no summation) we single out contributions  with the three
distinguished marked points belonging to the same joint of the curve.
The sub-sum turns out to coincide with $\Delta_{\a}^{1/2}$.
We refer to \cite{G1,G2} for further details of this 
{\em materialization} phenomenon in the theory of canonical coordinates. 
Our computation of higher genus potentials  
uses some of these results along with the standard fixed point localization
technique \cite{Kn1, GP} in the moduli spaces of stable maps. 
%In the discussion below we
%assume that the reader has some experience of working with localization 
%formulas in spaces of stable maps (see, for instance, \cite{Kn1,GP}).
   
The edge factors mentioned earlier are identified with $V_{kl}^{ij}$. First, 
in the fixed point expression for 
$e_{i}\lan \frac{\phi_{i}}{\x-c},\frac{\phi_{i}}{z-c}\ran$  
we single out contributions of those fixed point where the first and the 
last marked points belong to the same joint of the curve.
The sum of such contributions turns out to coincide 
with $e^{u^{i}(1/\x+1/z)}/(z+\x)$ (see \cite{G2,G3}). 
Therefore this expression occurs 
in the localization formula for the one-point descendent 
$\lan \phi_{\a}, \phi_i/(z-c)\ran$ as the factor responsible for the
contributions of the joints carrying the last marked point. The variable $\x$
is to be replaced by the character of the torus action on the leg approaching
the joint from the direction of the first marked point. 
Thus the dependence of the descendent on $z$ is 
transparent from the expansion of the factor: 
$e^{u^i/z}[\sum e^{u^i/\x}(-z)^k/\x^{k+1}]$.    
We conclude that the matrix 
$[\ \lan \phi_{\a},\frac{\phi_j}{z-c}\ran \sqrt{e_j}\ ]$
(normalized this way) is a unitary solution $S$ of the part (b) of 
Proposition.
It is one of the solutions described by the part (b) of Proposition.
Among the total class of solutions (see part (c)), it is characterized by the 
property that the series $R(z)$ turns into $1$ in the limit of classical equivariant cohomology,
that is when contributions of all non-constant stable maps are neglected. Eventually we will
have to change this normalization of the solution $S$ in order 
to compensate the effect of Hodge integrals in localization formulas.

Processing similarly contributions of the joints carrying
the first and last marked points in 
localization formulas for the two-point descendent (\ref{9}),
we extract the edge factors mentioned above:
\begin{equation}  \label{11} 
\lan \frac{\phi_i}{z-c},\frac{\phi_j}{w-c}\ran \sqrt{e_ie_j} = 
e^{u^i/z+u^j/w}[ \frac{\d_{ij}}{z+w}+\sum (-z)^k(-w)^l \ 
\text{(edge factor)}_{kl}^{ij} ] .\end{equation}
Taking into account (\ref{6}) and (\ref{10}) we conclude that the edge factors
are identified with the coefficients $V^{ij}_{kl}$ corresponding to
the solution $S$. Note that the weights $e^{u^j/\x}$ are incorporated 
into the edge factors.  

Computing contributions of vertices, 
denote by $\x^i_r, \ r=1,...,\dim_{\CC} X$, 
the characters of the torus action on the tangent space to $X$ at the 
fixed point with the index $i$. The localization formulas require
the following intersection numbers in the Deligne-Mumford
spaces:
\begin{equation} \label{12} 
\sum_{n=0}^{\infty} \frac{e_i^{-1}}{n!}
\int_{\M_{g,m+n}} \frac{\prod_{s=1}^g\prod_r (\x_r^i-\gr_s)}
{(x_1-c^{(1)})...(x_m-c^{(m)})}\w Q(c^{(m+1)})\w...\w Q(c^{(m+n)}) . 
\end{equation}
Here $\gr_1,...,\gr_g$ are Chern roots of the 
{\em Hodge bundle} with the fiber
$H^1(\Sigma, {\mathcal O}_{\Sigma})^*$, and $x_1,...,x_m$ are formal variables.
In localization formulas, these variables are replaced by some $\x_r^i$
(or their fractions), the characters of the torus action on 
the $m$ edges adjacent to the vertex. The formula (\ref{1}) accounts for this
substitution by matching the factors $c^{k}x^{-k-1}$ in (\ref{12}) with
the corresponding edge factors $V_{k...}^{i...}$  in (\ref{11}).  

The series $Q(c)=Q^i_0+Q^i_1c+...$ is to be substituted
in the localization formulas by the localization factor of the tail  
approaching the $i$-th fixed point, and the next task is to interpret 
the factor in terms of abstract Frobenius structures. For this,
we notice that the same series $Q(c)$ occurs --- in the same role ---
in fixed point localization of genus $0$ invariants. In particular,
the one-point descendent 
$\lan \frac{\phi_i}{z-c} \ran  = z \lan 1, \frac{\phi_i}{z-c}\ran$
is written as
\[ \frac{z}{e_i}+\frac{Q(-z)-Q(0)}{e_i}+\sum_{n=2}^{\infty}\frac{1}{n!}
\int_{\M_{0,1+n}}
Q(c^{(1)})\w...\w Q(c^{(n)})\w \ev_{n+1}^*\frac{\phi_i}{z-c} .\]
When $Q^i_0=0$, it coincides with $(z+Q(-z))/e_i$. This can be achieved 
by moving the series $Q$ by the {\em string flow}, and the time needed in order
to make $Q^i_0=0$ 
is exactly $-u_i$ (see \cite{G1}, Section $12$, or \cite{G2}).
Furthermore, both the descendent and the potential (\ref{12}) are 
eigenfunctions of 
the string operator $\p/\p Q_0-\sum Q_{k+1}\p/\p Q_k$   
(with the eigenvalues $1/z$ and $1/x_1+...+1/x_m$) and of the {\em dilaton
operator} $\p/\p Q_1-\sum Q_k \p/\p Q_k$ 
(with the eigenvalues $-1$ and $2g-2+m$ respectively). Thus, moving along the
string flow during the time interval $-u^i$ and then along the
dilaton flow during the time interval $\ln \sqrt{\Delta_i}$ we make 
$Q^i_0$ and $Q^i_1$ vanish and find the final values $Q^i_k=T^i_k$ from
(\ref{7}). The toll to pay consists of the factor $\Delta_i^{g-1+m/2}$ 
distributed in (\ref{1}) among vertices and edges and the 
weights $\exp u^i/\x_r^i$ already incorporated, as we remarked earlier, 
into $V^{i...}_{k...}$. 

\subsection{ Compensating constants.}

Yet, with our current definition
of $V_{kl}^{ij}$ and $T_k^i$ the formula (\ref{1}) would represent correctly 
the fixed point localization of higher genus potentials only if the Hodge
factors in (\ref{12}) were replaced with the factor 
$\prod_{s,r}\x_r^i=e_i^g$ (which cancels with other occurrences 
of $e_i$ here and there).
The Hodge factors should be digested as follows. Let $N_k$ denote Newton
symmetric polynomials. It is known \cite{FP} that $N_{2k}(\gr)=0$. We rewrite
\[ \prod_{s,r} (\x_r^i-\gr_s)=e_{i}^{g} 
\exp [ - \sum_{k=1}^{\infty} N_{2k-1}(1/\x^i) N_{2k-1}(\gr)/(2k-1) ] .\] 
Let us redefine the fundamental solution 
$S=[\lan \phi_{\a},\phi_{i}/(z-c)\ran e_{i}^{1/2}]$ 
using the ambiguity
described in the part (c) of Proposition:
\[ \ ^{new}S_{\a}^i:= S_{\a}^i 
\exp [-\sum z^{2k-1}\frac{N_{2k-1}(1/\x^i)}{2k-1}\frac{B_{2k}}{2k}].\]
Here $B_{2k}$ denote Bernoulli numbers, 
$z/(\exp z-1)=1-z/2+\sum_k z^{2k}B_{2k}/(2k)!$
The coefficients $V^{ij}_{kl}$ in (\ref{6}) are redefined accordingly.

\medskip

{\bf Theorem 1.} 
{\em In equivariant Gromov -- Witten theory for Hamiltonian tori
actions with isolated fixed points, we obtain the higher genus potential
in the form (\ref{1}) by taking $\ ^{new}S$ on the role 
of the fundamental solution $S$ in (\ref{6}) and (\ref{7}).} 

\medskip

{\em Remark.} According to the Proposition, the unitary solution $\ ^{new}S$ is 
charachterized by the condition that the corresponding series $\ ^{new}R(z)$ turns into
the diagonal matrix of the compensating constants  
$\exp [-\sum z^{2k-1}\frac{N_{2k-1}(1/\x^i)}{2k-1}\frac{B_{2k}}{2k}] $
in the limit of {\em classical} equivariant cohomology. Thus Theorem $1$, under its hypotheses,
coincides with Conjecture $1$ where the solution $S$ defined on the basis of Proposition
is normalized in this particular way.   
 
\medskip

{\em Example 6.} In genus $1$, the differential of the GW-potential
was computed by fixed point localization in \cite{G2}. In our current
notation 
\[ dF_X^1= \sum_{i}[\ \frac{V_{00}^{ii}}{2}du^i
-\frac{N_1(1/\x^i)}{24}du^i+\frac{d\Delta_i}{48\Delta_i}\ ] .\]
The first summand represents contributions of cycles of rational curves,
that is of graphs with one vertex (of type $(g,m)=(0,3)$) and one edge. 
The other two summands come from (\ref{12}) with $(g,m)=(1,1)$.
The middle term is due to the Hodge integral 
$\int_{\M_{1,1}} \sum \gr_s =1/24$. It can be interpreted as contributions
of cycles of rational curves shrinking to a point and is incorporated into
the first term as $\ ^{new}V_{00}^{ii}=V_{00}^{ii}-N_1(1/\x^i)/12$. 
This change of notation agrees with the theorem since $B_2/2=1/12$.
We arrive at the conjecture \cite{G2} making sense for arbitrary semisimple
Frobenius manifolds:
\[ dF^1=\sum_{i}[\ \frac{V_{00}^{ii}}{2}du^i + 
\frac{d\Delta_i}{48\Delta_i}\ ]. \]
In the case of conformal Frobenius structures the conjecture was proved in
\cite{DZ1} by showing that this is the only homogeneous 
formula that agrees with Getzler's equation \cite{Ge}.     

\subsection{  Hodge intersection numbers.} 

We have already explained why the formula
(\ref{1}) for higher genus potentials would arise if the Hodge factors
in the vertex contributions (\ref{12}) were neglected. To derive the theorem
it remains to prove that the effect of Hodge factors is correctly 
accounted by the modification $S \mapsto \ ^{new}S$. 
For this, let us introduce the
generating function for Hodge intersection numbers:
\begin{equation} \label{13} 
\gl (\h; Q; s_1,s_2,...)= \exp \{ \sum_{g=0}^{\infty} 
\h^{g-1} {\mathcal H }^{g,n}_{pt}(Q, s_1,s_2,...) \} \end{equation}
where
\[ {\mathcal H}^{g,n}_{pt}:=
 \sum_{n=0}^{\infty} \frac{1}{n!}\int_{\M_{g,n}} 
Q(c^{(1)})\w ... \w Q(c^{(n)})\w e^{\sum s_k N_{2k-1}(\gr)/(2k-1)!} . \]
We can introduce a family of fake higher genus potentials depending on
the parameters $\{ s_k^i\}$ by replacing the factors 
$\gt (\h \Delta_i; Q^i)$ in (\ref{1}) with
$\gl (\h \Delta_i; Q^i; s^i_1,s^i_2,...)$. The actual higher genus
potential corresponds to $s^i_k=-(2k-2)!\ N_{2k-1}(1/\x^i)$. 
We claim that {\em the $s$-parametric deformation of (\ref{1}) is 
identified with the $a$-parametric deformation of the fundamental 
solution $S$ described in the part (c) of Proposition by taking
$a_k^i=B_{2k} s_k^i/(2k)!$} This obviously implies the Theorem.

Following (\ref{6}) and (\ref{7}) with
scalar $R(z)=\exp (a_1z+a_2z^3+...) $ and $\Delta=1$, we introduce 
the operator  
$P(a_1,a_2,...)=\frac{1}{2}\sum v_{kl}(a)\p_{\tilde{Q}_k}\p_{\tilde{Q}_l}$
where
\begin{equation} \label{14}
 \frac{1}{z+w}+\sum v_{kl}(a_1,a_2,...)\ (-z)^k(-w)^l:=
\frac{\exp \{\sum 
a_k (z^{2k-1}+w^{2k-1})\} }{z+w} \end{equation}
and define a substitution $\tilde{Q}(Q,s)$ by
\begin{equation} \label{15} z+\tilde{Q}(-z):=[z+Q(-z)]
\exp [\sum  a_k z^{2k-1}]. \end{equation}

\medskip

{\bf Lemma.}  
\[ \gl (\h; Q; s_1,s_2,...) = 
[e^{\h P ( \frac{B_2}{2!}s_1,\frac{B_4}{4!}s_2,...)} 
\gt (\h ; \tilde{Q}) ]  
_{\tilde{Q}=\tilde{Q}(Q, \frac{B_2}{2!} s_1, \frac{B_4}{4!} s_2, ...)} \]

\medskip

Our claim follows formally from Lemma. Indeed, the $a$-parametric
modification $\ ^{new}S=S\exp (\sum a_kz^{2k-1})$ of the fundamental solution 
affects the values $Q_k^i=T_k^i$ by some linear transformation and also 
changes coefficients $V_{kl}^{ij}$ of the differential operator in the exponent
of (\ref{1}). Instead of changing the values $T_k^i$ one can make the
change of the variables $Q^i\mapsto \tilde{Q}^i$ and leave the values 
$Q_k^i=T_k^i$ unchanged. The change of variables coincides with (\ref{15}). 
The same change of variables in the differential operator accounts 
for the most of the change in the coefficients $V_{kl}^{ij}$. The only
remaining discrepancy comes from the term $\d_{ij}/(z+w)$ in (\ref{6}) 
and is determined by (\ref{14}) as
$\ ^{new}V_{kl}^{ij}=V_{kl}^{ij}+\d_{ij} v_{kl}(a_1^i,a_2^i,...)$. Thus
$\sum_i \h \Delta_i P(a^i) $ is added to the differential operator 
in the exponent of (\ref{1}). According to Lemma the modification
is equivalent to using $\gl (\h; Q; s) $'s (instead of $\gt (\h ; Q)$ in
(\ref{1})) when $a_k^i=B_{2k}s_k^i/(2k)!$  

\medskip

{\em Proof of the lemma.}
It is known \cite{F}, at least in principle, how to compute 
$\gl $ in terms of $\gt$ using Mumford's Grothendieck - Riemann - Roch 
formula \cite{Mu} for the
Chern character $-\sum N_{2k-1}(\gr)/(2k-1)!$ of the Hodge bundle. 
Moreover, the formula is interpreted in \cite{FP} as the PDE-system
\begin{equation} \label{16} 
\frac{\p}{\p s_m}\gl = \frac{B_{2m}}{(2m)!} (\h D_m + L_m) \gl,\ 
m=1,2,... \end{equation} 
where $ D_m:= \frac{1}{2}\sum_{k+l=2m-2}(-1)^k \p_{Q_k}\p_{Q_l}$
and $L_m:=\p_{Q_{2m}}-\sum_{k=0}^{\infty} Q_k \p_{Q_{k+2m-1}}$.
The operators $\h D_m+L_m$ commute pairwise. The vector fields $L_m$ on
the space of power series $Q(c)=Q_0+Q_1c+...$ are linear with
respect to the origin shifted to $c$. In fact they are given by the 
operators of multiplication by $-c^{2m-1}$. Therefore $L_m$
commute themselves and define the flow (\ref{15}).  
Furthermore, for functions $f(\tilde{Q})$ we find by differentiation
that  $[\frac{\p}{\p a_m} (P f)] 
(\tilde{Q}(Q,a))= D_m [ f(\tilde{Q}(Q,a) ]$. The lemma follows: both 
sides satisfy the same PDE system (\ref{16}) and coincide at $s=0$.

\section{Generalization to gravitational descendents}
 
The genus $g$ descendent GW-potential
of $X$ is a formal function on the {\em space of curves} 
$\t=t_0+t_1c+t_2c^2+...$ in $H$ defined by
\begin{equation} \label{17}
{\mathcal F}^g_X (\t):=\sum_{n,d}\frac{q^d}{n!}\int_{[X_{g,n,d}]}
\ev_1^*\t (c^{(1)})\w ... \w \ev_n^*\t(c^{(n)}) . \end{equation}
Here $c^{(i)}$ is the $1$-st Chern class of the universal cotangent line
over $X_{g,n,d}$ at the $i$-th marked point, and $\ev_i^*$ acts on coefficients
$t_m$ of the series $\t$. We present here a conjectural formula for 
higher genus descendent potentials that makes sense for arbitrary 
semisimple Frobenius structures.
For this, we have to review the
construction \cite{Db} of genus $0$ descendents of Frobenius manifolds.  

\subsection{  Descendents in genus 0.}

One starts with a fundamental solution $1+z^{-1}S_1+z^{-2}S_2+...$ to the system (\ref{3})
satisfying the unitary condition and takes it on the role of the one-point descendent (\ref{8}):
\[ (\lan \phi_{\a}, \phi_{\mu}/(z-c)\ran g^{\mu\b}) := 
1+ \sum_{k>0} z^{-k} (\lan \phi_{\a}, \phi_{\mu} c^{k-1}\ran ' g^{\mu\b}) := 
1 + \sum_{k>0} z^{-k}S_k .\]
We emphasize that the $1/z$-series solution is considered disjoint from 
the asymptotical solution $S=\Psi R(z) \exp (U/z)$ of the Proposition.
In particular, in equivariant GW-theory the one-point descendent correlators, 
defined intrinsicly, form the fundamental solution series in question, and this definition 
is not affected by the modification $R\mapsto \ ^{new}R$ of integration constants in the
series $R$. 
           
Next, one introduces the $2$-point descendent (\ref{9}) using (\ref{10}): 
\begin{equation} \label{18} 
\lan \frac{\phi_{\a}}{z-c},\frac{\phi_{\b}}{w-c}\ran = \frac{g_{\a\b}}{z+w}+
\sum \frac{\lan \phi_{\a}c^m,\phi_{\b}c^l\ran '}{z^{m+1}w^{l+1}}
:= \lan \phi_{\mu},\frac{\phi_{\a}}{z-c}\ran \frac{g^{\mu\nu}}{z+w}
\lan \phi_{\nu},\frac{\phi_{\b}}{w-c}\ran  .\end{equation} 
The singular term is present to make the sum satisfy the string equation
but it makes the symbol $\lan \cdot,\cdot \ran$ not entirely bilinear. We
use here the notation $\lan \cdot,\cdot\ran '$ for the honest bilinear 
$2$-point descendents.
 
Furthermore, one considers the map 
\begin{equation} \label{19} \t\mapsto t(\t)=crit\ \lan \t(c)-c, 1\ran (t) 
\end{equation}   
from the curve space to the Frobenius manifold defined by taking the critical 
point of the function $\lan \t(c)-c, 1\ran:= (t_0,t)+\lan \t(c)-c,1\ran '$
of $t\in H$ depending linearly
on the parameter $\t=t_0+t_1c+...$. One can show that the equation of the 
critical point takes on the form 
$t^{\a}=t_0^{\a}+g^{\a\mu}\lan \phi_{\mu}, (\t(c)-\t(0))/c\ran (t) $
and thus admits a unique formal solution which turns into $t=t_0$ when 
$t_1=t_2=...=0$. Finally one puts 
\begin{equation} \label{20}
{\mathcal F}^{0}(\t)=\frac{1}{2}\lan \t(c)-c,\t(c)-c\ran '(t(\t)) 
\end{equation}
As it is shown in \cite{Db}, the formula (\ref{20})  
agrees with the {\em string equation} and the genus $0$ {\em topological
recursion relation} and is the only deformation of 
${\mathcal F}^0|_{t_1=t_2=...=0}=F^0(t_0)$ satisfying these conditions.  
Also, (\ref{20}) agrees with the {\em dilaton equation} and is consistent 
with the definition (\ref{18}):
\begin{equation} \label{21}
\p_{t^{\a}_m}{\mathcal F}^0(\t)=\lan \phi_{\a}c^m,\t(c)-c\ran '(t(\t)),\ 
\p_{t^{\a}_m}\p_{t^{\b}_l}{\mathcal F}^0(\t)=
\lan \phi_{\a}c^m,\phi_{\b}c^l\ran ' (t(\t)).\end{equation}
We emphasize that all the two-point descendent correlators 
$\p_{t^{\a}_m}\p_{t^{\b}_l}{\mathcal F}^0$ depend on the infinitely many variables 
$\t$ only via the substitution $t(\t)$.

\subsection{Descendents in higher genus.} 

Our proposal for higher genus descendent
potential has the same form as (\ref{1}):
\begin{equation} \label{22} 
e^{\sum_{g\geq 2}\h^{g-1}{\mathcal F}^{g}(\t)}=
[\ e^{\frac{\h}{2}\sum \V_{kl}^{ij}\sqrt{\Dl_i\Dl_j}\p_{Q^i_k}\p_{Q^j_l}}
\prod_j \gt (\h \Dl_j; Q^j)\ ]_{Q^i_k=\T^i_k} \ .\end{equation}
The functions $\V_{kl}^{ij},\Dl_i,\T^i_k$ on the curve space are defined 
near a semisimple point $\t(0)$ in terms of genus $0$ descendents. 
The definitions are motivated by computation of higher genus descendent 
potentials in equivariant GW-theory in the presence of the torus acting on the
target space with isolated fixed points. The result of the computation is stated 
in Theorem $2$ below. The proof of Theorem $2$ follows is identical to
the proof of Theorem $1$ given in Section $2$, with one deviation which is discussed 
presently. In the part of the text preceeding the formulation of Theorem $2$
we assume that the reader is familiar with the details of Section $2$.

Let us remind from Section $2$ that the formula (\ref{22}) with the combinatorial structure 
of a graph sum originates from the technique of fixed point localization in moduli spaces
of stable maps, and that the edge and tail factors are to be extracted from the genus $0$
descendents.   
In particular the edge factors in localization formulas are
extracted from the expansion (\ref{11}) for $2$-point correlators on the
curve space. Due to (\ref{21}) all such $2$-point correlators coincide
with the corresponding $2$-point descendents on $H$ lifted to the
curve space by the change of variables (\ref{19}). 
This also applies to $u^i=u^i(t(\t))$ (which can be described 
via $2$-point correlators, see \cite{G1,G2}) and therefore --- to the edge factors: 
\begin{equation} \label{23} \V_{kl}^{ij}(\t)=V_{kl}^{ij}(t(\t)) \ 
\text{where $t(\t)$ is defined by (\ref{19})}.\ \end{equation} 
We stress that in the present context of fixed point localization the edge factors 
$V_{kl}^{ij}$ will eventually have to be the same as in Theorem $1$ 
(\this based on the solution $\ ^{new}S$ modified by the Bernoulli constants)  
in order to compensate the effect of Hodge integrals.

Similarly, the functions $\Dl_i$ and $\T^i_k$ in the localization formulas
are found from the expansion of the $1$-point correlator on the curve space:
\[ \sqrt{e_i}\ \lan \frac{\phi_i}{z-c}\ran (\t)
=\frac{e^{u^i/z}}{\sqrt{\Dl_i}}(z+\sum \T_k^i(-z)^k) .\]
However --- and this is the point that makes the difference --- 
$\lan \frac{\phi_i}{z-c}\ran $ no longer coincides with 
the $2$-point correlator $z\lan 1, \frac{\phi_i}{z-c}\ran $ and respectively $\Dl_i$ and
$\T^i_k$ are not obtained from $\Delta_i$ and $T^i_k$ by the substitution $t(\t)$. 
We need (\ref{18}--\ref{21}) in order to interpret them in terms of 
abstract Frobenius structures. Let us recall that the asymptotical solution $S^i_{\mu}(z)$
in the context of fixed point localization actually coincides with 
$\lan \phi_{\mu},\frac{\phi_{i}}{z-c}\ran \sqrt{e_i}$. We have  
\[ \sqrt{e_i}\lan \frac{\phi_i}{z-c}\ran = 
\ \sum S^i_{\mu}(z) g^{\mu\nu}\oint \ \lan \phi_{\nu},\frac{\phi_{\a}}{w-c}\ran \ \frac
{(\t^{\a}(w)-\d^{\a}w)}{2\pi i (z+w)}\ dw .\] 
Computing the countur integral we arrive at the formula
\begin{equation} \label{24} \begin{array}{l} 
\frac{e^{u^i/z}}{\sqrt{\Dl_i}}(z+\sum\T^i_k(-z)^k)=  
 \ \sum_{\mu\nu} S^i_{\mu}(z)_{t=t(\t)} g^{\mu\nu} \times 
 \{ \lan \phi_{\nu},1,\t(c)-c\ran + \\ \\ (-z)\lan \phi_{\nu},1,1,\t(c)-c\ran+
(-z)^2\lan\phi_{\nu},1,1,1,\t(c)-c\ran+...\}_{t=t(\t)} .
\end{array} \end{equation}    
Here $\lan \phi_{\nu},1,...,1, f(c)\ran (t)$ coincide with 
multiple $t$-derivatives of $\lan \phi_{\nu}, f(c)\ran (t)$ in the 
direction of the vector $1$. Here the 
notation $S^i_{\mu}(z)=e^{u^i/z}(\sum (R_k)_j^iz^k)\Psi^j_{\mu}$ in the context of
localization formulas should eventually refer to $\ ^{new}S$, 
the fundamental solution matrix modified by the Bernoulli constants.  

{\em We take (\ref{23}) and (\ref{24}) on the role of definitions for 
$\V_{kl}^{ij}, \T^i_k$ and $\Dl_i$ in the formula (\ref{22}).} 

By definition $\T^i_1=0$ while $\T^i_0=0$ follows from the criticality 
condition in (\ref{19}). It is straightforward to check that the 
definition reduces to (\ref{6},\ref{7}) when $t_1=t_2=...=0$. With these
definitions in force, and the compensating constants in $R(z)$ in place, 
we arrive at the following theorem.
 
\medskip

{\bf Theorem 2.} {\em The formula (\ref{22}) for higher genus descendent 
potentials holds true in equivariant GW-theory of Hamiltonian tori actions 
with isolated fixed points.}

\medskip

{\em Example 7.} We compute from (\ref{24}) that
\[ \Dl_i^{-1/2} (\t)=\Delta_i^{-1/2} [\sum \frac{\p u^i}{\p t^{\mu}}
g^{\mu\nu} \lan \phi_{\nu}, 1, 1, c-\t (c)\ran ] (t(\t)) .\]
Along the lines of Example $6$ we get $ d{\mathcal F}^1=$
\[ =\sum (\frac{\V_{00}^{ii}}{2}du_i+\frac{d\Dl_i}{48\Dl_i})=
d \{ F^1(t(\t)) - \frac{1}{24} 
\ln [ \prod_i \frac{\p u^i}{\p t^{\mu}}g^{\mu\nu}
\lan \phi_{\nu},1,1,c-\t(c)\ran ] \} .\]
This answer actually coincides with the well-known result \cite{DW}
\[ {\mathcal F}^1(\t)=F^1(t(\t))+\frac{1}{24}
\ln \det [ \frac{\p t^{\mu}}{\p t^{\nu}_0} ] .\]  
Indeed, differentiating the criticality condition 
$\lan \phi_{\d},1,c-\t(c)\ran=0$ in (\ref{19}) we find
that $g^{\a\gel}\lan \phi_{\gel},\phi_{\b},1,c-\t(c)\ran (t(\t))$ form the
matrix inverse to $[\p t^{\mu}/\p t_0^{\nu}]$. On the other hand,
the genus $0$ topological recursion relation (or WDVV-equation) implies
\[ g^{\a\gel}\lan \phi_{\gel},\phi_{\b},1,c-\t(c)\ran =
g^{\a\gel}\lan \phi_{\gel},\phi_{\b},\phi_{\mu}\ran g^{\mu\nu}
\lan \phi_{\nu}, 1, 1, c-\t(c)\ran .\]
In other words, $[\p t/\p t_0]^{-1}$ coincides with the linear combination
with coefficients $g^{\mu\nu}\lan \phi_{\nu},1,1,c-\t(c)\ran$ of the commuting
matrices $[ g^{\a\gel}F^0_{\gel\mu\b} ]$  of quantum multiplication  
operators $\phi_{\mu} \bullet$. Thus the eigenvalues of the matrix 
are the linear combinations $(\p u^i/\p t^{\mu}) g^{\mu\nu} \lan \phi_{\nu},
1,1,c-\t(c)\ran $, and the determinant is their product.

\medskip
       
{\em Example 8.} Consider GW-theory with the target space $X=pt$.
Then $u=t$, $\lan 1, 1/(z-c)\ran=\exp (t/z)$ and respectively $\Delta =1$ and 
$V_{kl}=0$. We find the RHS in (\ref{22}) equal to $\gt (\h \Dl ; \T )$ 
with $\Dl $ and $\T_k$ computed as follows. We have
$f(t; \t):=\lan 1,1,\t(c)-c\ran (t)=\sum t_kt^k/k!-t$. The relation 
(\ref{19}) turns into $f(t(\t);\t)=0$, while $\Dl^{-1/2}=-f'(t(\t);\t)$, and
\[ \T_1-1=f'(t(\t);\t) \sqrt{\Dl}, \ \T_2=f''(t(\t),\t)\sqrt{\Dl},\ 
\T_3=f'''(t(\t);\t)\sqrt{\Dl},...\]
Note that $(t_k-\d_{k,1})\mapsto \p^k f(t;\t)/\p t^k$ is the string flow on the
curve space so that $\T$ is obtained from $\t $ by applying the
string flow until $t_0=0$ and then applying the dilaton flow until $t_1=0$.
The potentials ${\mathcal F}^g_{pt}(\t)$ with $g=0,1$ vanish when 
$t_0=t_1=0$, and for $g\geq 2$ are preserved by the string flow and are
homogeneous of degree $2-2g$ with respect to the dilaton flow. 
We conclude that indeed $\gt (\h \Dl ; \T )$ coincides with 
$\exp \sum_{g\geq2} \h^{g-1} {\mathcal F}^g_{pt}(\t)$.

Finally, the formulas (\ref{22}--\ref{24}) agree with our
Lemma about Hodge intersection numbers in the following sense: 
the Lemma follows formally from our claim that, in the current
setting with descendents as well, the $s$-deformation (\ref{13}) 
of (\ref{2}) is compensated by the modification 
$\exp (u/z) \mapsto \exp (u/z + B_2s_1z/2! + B_4 s_2z^3/4!+...)$
described in the part (c) of Proposition.

\medskip

{\bf Conjecture $2$.} {\em With the notations (\ref{2},\ref{23},\ref{24}) in force,
the formula (\ref{22}) represents higher genus gravitational descendents of compact
symplectic manifolds with generically semisimple quantum cup-product.}

\medskip

{\em Remark.} Generally speaking, the task of deriving Conjectures $1$ and $2$ 
from Theorems $1$ and $2$ by passing to the non-equivariant limit is open and non-trivial. 
We will show in \cite{G4} how to do this in the case of complex projective spaces and their 
products.

\subsection{Concluding remarks.}

The proposal (\ref{1},\ref{22}) should be exposed to further tests. 
We expect it to be consistent with any 
relations in cohomology of Deligne--Mumford spaces (see \cite{Ge,Ge2}
for some such ralations found in genus $0$ and $2$). 
In fact we hope that our Theorems $1$ and $2$
on equivariant GW-potentials impose
some constraints on topology of Deligne--Mumford spaces so tight that 
the corresponding results in abstract semisimple GW-theory 
would follow. Respectively, it should be interesting to make such 
constraints as explicit as possible and  
to understand better the geometrical structure on Frobenius manifolds 
encrypted by (\ref{1}) and (\ref{22}). To this end, we should say that
the formulas (\ref{1},\ref{22}) can be rewritten differently. Using the 
Fourier transform they can be given a form of path integrals. Substituting
matrix Airy integrals for the Kontsevich -- Witten function (\ref{2}) we
can relate the formulas to multi-matrix models.
Perhaps, the most useful is the representation-theoretic formulation \cite{G4},
which automatically restores the (descendent) potentials of genus $0$ and $1$ 
and yields a formula for the complete tau-function 
$\exp [\sum_{g\geq 0}h^{g-1}{\mathcal F}^g]$. It enables us to show 
(see \cite{G4}) that the proposal agrees with the Virasoro constraints 
\cite{DZ2} and to derive the ``Virasoro conjecture'' for GW-invariants of 
complex projective spaces and their products.
We also expect the formula for the tau-function to be helpful in the 
construction of the bihamiltonian structure of the KdV-like integrable 
hierarchy whose approximations are studied in \cite{Db,DZ1,DZ2,EGX}.

\section*{Acknowledgments}  
  
The author would like to thank the National Science Foundation 
for financial support. 
A major part of the text has been completed during our stay at
Ecol\'e Polytechnique and then at Schr\"odinger Institute whose hospitality
is thus gratefully acknowledged. 
Many thanks are due to C. Faber and R. Pandharipande
for their interest and especially for their generous help with some intersection 
numbers on Deligne -- Mumford spaces needed in the early stages of this work to 
uncover the relation between the compensating constants and Bernoulli numbers.

\enddocument